\documentclass[12pt]{amsart}
\usepackage{amssymb}
\usepackage{verbatim}
\usepackage{hyperref}
\usepackage{color}

\begin{document}

\newtheorem{theorem}{Theorem}    
\newtheorem{proposition}[theorem]{Proposition}
\newtheorem{conjecture}[theorem]{Conjecture}
\def\theconjecture{\unskip}
\newtheorem{corollary}[theorem]{Corollary}
\newtheorem{lemma}[theorem]{Lemma}
\newtheorem{sublemma}[theorem]{Sublemma}
\newtheorem{observation}[theorem]{Observation}
\theoremstyle{definition}
\newtheorem{definition}{Definition}
\newtheorem{notation}[definition]{Notation}
\newtheorem{remark}[definition]{Remark}
\newtheorem{question}[definition]{Question}
\newtheorem{questions}[definition]{Questions}
\newtheorem{example}[definition]{Example}
\newtheorem{problem}[definition]{Problem}
\newtheorem{exercise}[definition]{Exercise}

\numberwithin{theorem}{section} \numberwithin{definition}{section}
\numberwithin{equation}{section}

\def\earrow{{\mathbf e}}
\def\rarrow{{\mathbf r}}
\def\uarrow{{\mathbf u}}
\def\varrow{{\mathbf V}}
\def\tpar{T_{\rm par}}
\def\apar{A_{\rm par}}

\def\reals{{\mathbb R}}
\def\torus{{\mathbb T}}
\def\heis{{\mathbb H}}
\def\integers{{\mathbb Z}}
\def\naturals{{\mathbb N}}
\def\complex{{\mathbb C}\/}
\def\distance{\operatorname{distance}\,}
\def\support{\operatorname{support}\,}
\def\dist{\operatorname{dist}\,}
\def\Span{\operatorname{span}\,}
\def\degree{\operatorname{degree}\,}
\def\kernel{\operatorname{kernel}\,}
\def\dim{\operatorname{dim}\,}
\def\codim{\operatorname{codim}}
\def\trace{\operatorname{trace\,}}
\def\Span{\operatorname{span}\,}
\def\dimension{\operatorname{dimension}\,}
\def\codimension{\operatorname{codimension}\,}
\def\nullspace{\scriptk}
\def\kernel{\operatorname{Ker}}
\def\ZZ{ {\mathbb Z} }
\def\p{\partial}
\def\rp{{ ^{-1} }}
\def\Re{\operatorname{Re\,} }
\def\Im{\operatorname{Im\,} }
\def\ov{\overline}
\def\eps{\varepsilon}
\def\lt{L^2}
\def\diver{\operatorname{div}}
\def\curl{\operatorname{curl}}
\def\etta{\eta}
\newcommand{\norm}[1]{ \|  #1 \|}
\def\expect{\mathbb E}
\def\bull{$\bullet$\ }
\def\C{\mathbb{C}}
\def\R{\mathbb{R}}
\def\Rn{{\mathbb{R}^n}}
\def\Sn{{{S}^{n-1}}}
\def\M{\mathbb{M}}
\def\N{\mathbb{N}}
\def\Q{{\mathbb{Q}}}
\def\Z{\mathbb{Z}}
\def\F{\mathcal{F}}
\def\L{\mathcal{L}}
\def\S{\mathcal{S}}
\def\supp{\operatorname{supp}}
\def\dist{\operatorname{dist}}
\def\essi{\operatornamewithlimits{ess\,inf}}
\def\esss{\operatornamewithlimits{ess\,sup}}
\def\xone{x_1}
\def\xtwo{x_2}
\def\xq{x_2+x_1^2}
\newcommand{\abr}[1]{ \langle  #1 \rangle}

\newcommand{\Norm}[1]{ \left\|  #1 \right\| }
\newcommand{\set}[1]{ \left\{ #1 \right\} }
\def\one{\mathbf 1}
\def\whole{\mathbf V}
\newcommand{\modulo}[2]{[#1]_{#2}}

\def\scriptf{{\mathcal F}}
\def\scriptg{{\mathcal G}}
\def\scriptm{{\mathcal M}}
\def\scriptb{{\mathcal B}}
\def\scriptc{{\mathcal C}}
\def\scriptt{{\mathcal T}}
\def\scripti{{\mathcal I}}
\def\scripte{{\mathcal E}}
\def\scriptv{{\mathcal V}}
\def\scriptw{{\mathcal W}}
\def\scriptu{{\mathcal U}}
\def\scriptS{{\mathcal S}}
\def\scripta{{\mathcal A}}
\def\scriptr{{\mathcal R}}
\def\scripto{{\mathcal O}}
\def\scripth{{\mathcal H}}
\def\scriptd{{\mathcal D}}
\def\scriptl{{\mathcal L}}
\def\scriptn{{\mathcal N}}
\def\scriptp{{\mathcal P}}
\def\scriptk{{\mathcal K}}
\def\frakv{{\mathfrak V}}

\title[Sharp Weighted Bounds for Multilinear fractional Maximal type Operators]
{Sharp Weighted Bounds for Multilinear fractional Maximal type Operators with Rough Kernels}

\author[Ting Mei]{Ting Mei} 

\author[Qingying Xue]{Qingying Xue}

\author[Senhua Lan]{Senhua Lan}
%
\thanks{ Corresponding author: Qingying Xue.
This work was partially supported by
 NSF of China (Grant No. 10931001), Beijing Natural Science Foundation
(Grant: 1102023). The third named author was supported partly by
NSF of China (Grant No.11171137) and NSF of Zhejiang Province (Grant No. LY12A01011).}
\address{Ting Mei\\School of Mathematical
Sciences\\Beijing Normal University\\Laboratory of Mathematics and
Complex Systems\\Ministry of Education\\Beijing, 100875\\
P. R. China}
\email{meiting1989@sina.cn}

\address{Qingying Xue\\School of Mathematical
Sciences\\Beijing Normal University\\Laboratory of Mathematics and
Complex Systems\\Ministry of Education\\Beijing, 100875\\
P. R. China
}
\email{qyxue@bnu.edu.cn}
\address{Senhua Lan\\Department of Mathematics\\Lishui University \\Li Shui, 323000\\
P. R. China
}
\email{senhualan@sina.com}
\maketitle

\begin{abstract}

In this paper, we will give the weighted bounds for multilinear fractional maximal type operators $\mathcal{M}_{\Omega,\alpha}$ with rough homogeneous kernels. We obtain a mixed $A_{(\vec{P},q)}-A_\infty$ bound and a $A_{\vec{P}}$ type estimate for $\mathcal{M}_{\Omega,\alpha}$. As an application, we give an almost sharp estimate for the multilinear fractional integral operator with rough kernels $\mathcal{I}_{\Omega,\alpha}$.

\end{abstract}

\section{Introduction} 

Multilinear Calder$\acute{o}$n-Zygmund operators were introduced and
first studied by Coifman and Meyer \cite{coif1}, \cite{coif2},
\cite{coif3}, and later on by Grafakos and Torres \cite{GT},
\cite{GT2}. The multilinear operators are natural generalizations of linear case. In recent years, the theory on multilinear
Calder$\acute{o}$n-Zygmund operators and related operators, such as multilinear singular integral, maximal and fractional maximal type operators, fractional integrals, have attracted much attentions. We begin by introducing the definition of the multilinear fractional maximal type operator as follows.
\begin{definition} \label{def 1.1}(\cite{LOPTT} or \cite{CX})
Given $\vec{f}=(f_1,\cdots,f_m)$, suppose each $f_i$ $(i=1,\cdots ,m)$ is locally integrable on $\R^n$. Then for any $x\in \R^n$, we define the multilinear fractional type maximal operator $\mathcal{M}_\alpha$ and the multilinear fractional integral operator $\mathcal{I}_{\alpha}$ by
\begin{align}\label{e 1.1}
\mathcal{M}_\alpha(\vec{f})(x)=\sup_{Q\ni x} \prod_{i=1}^m \frac{1}{|Q|^{1-\frac{\alpha}{mn}}} \int_Q |f_i(y_i)| \;dy_i,   \quad\quad for \ \ 0\le \alpha<n\end{align}
and
\begin{equation}\mathcal{I}_{\alpha}(\vec{f})(x) = \int_{{(\mathbb{R}^n)}^m}\frac{1}{{|(y_1, \cdots, y_m)|}^{mn - \alpha}}
\prod_{i = 1}^{m}f_i(x - y_i)\,d\vec{y},\quad\quad for \ \ 0< \alpha<n,\end{equation}where
the supremum in (\ref{e 1.1}) is taken over all cubes $Q$ containing $x$ in $\R^n$ with the sides parallel to the axes and $d\vec{y} = dy_1\cdots dy_m$.
\end{definition}

If $m=1$ in Definition \ref{def 1.1}, we simply denote $\mathcal{M}_\alpha$ by ${M}_\alpha$, $\mathcal{I}_{\alpha}$ by ${I}_{\alpha}$. It is well known that ${M}_\alpha$ and ${I}_{\alpha}$ play important roles in Harmonic Analysis and other fields, such as PDE, see e.g. \cite {riesz potential}, \cite{Fefferman good lambda}, \cite{multi frac grafakos}, \cite{MW}. A class of $A_{p,q}$ weights was first defined and weighted estimates of ${M}_\alpha$ and ${I}_{\alpha}$ were
considered by Muckenhoupt and Wheeden \cite{MW} in 1974. They proved that the fractional maximal operator $M_\alpha$ and the fractional integral operator ${I}_{\alpha}$ were of weak type $(L^1(\omega),L^{\frac{n}{n-\alpha},\infty}(\omega^{\frac{n}{n-\alpha}}))$ and of strong type $(L^p(\omega^p),L^q(\omega^q))$ if $p>1$, with$\frac 1q=\frac 1p-\frac \alpha n$ and $\alpha>0$. In 2010, if $w\in A_{p,q}$, Lacey et al \cite {{LMPT}} gave the sharp constant for $M_\alpha$ in the full range of exponents:\begin{equation}\|
{M}_\alpha\|_{L^q(w^q)}\le C[w]_{A_{p,q}}^{\frac{p'}{q}(1-\frac\alpha n)}\|f\|_{L^p(w^p)}, \quad\quad for \ \ 0\le \alpha<n, 1<p<n/\alpha.\label{1.3}\end{equation}
They pointed out that the exponent $\frac{p'}{q}(1-\frac\alpha n)$ is sharp and can not be improved. We noted that, for $\alpha=0$, the above result was first obtained by Buckley \cite{SB} for any $1<p<\infty$ and $\omega\in A_p$. Lacey et al \cite{LMPT} also obtained the following inequality for $I_\alpha$,
\begin{equation}\|
{I}_\alpha\|_{L^q(w^q)}\le C[w]_{A_{p,q}}^{\max\{1,\frac{p'}{q}\}(1-\frac\alpha n)}\|f\|_{L^p(w^p)}, \quad\quad for \ \ 0< \alpha<n, 1<p<n/\alpha.\label{1.4}\end{equation}Furthermore this estimate is sharp.

If $\alpha=0$ in (\ref{e 1.1}), then the operator $\mathcal M_0$ coincides with the new maximal operator $\mathcal{M}(\vec{f})$ defined by  Lerner, Ombrosi, P\'{e}rez, Torres and Trujillo-Gonz\'{a}lez in \cite{LOPTT}.
The authors in \cite{LOPTT} introduced the  following definition of multiple $A_{\vec{P}}$ weights.
\begin{definition}\label{def 1.2}(\cite{LOPTT})
Let $\vec{P}=(p_1,\cdots ,p_m), $ $\frac{1}{p}=\frac{1}{p_1}+\cdots+\frac{1}{p_m}$ with $1<p_1,\cdots ,p_m<\infty$, \, $\vec{\omega}=(\omega_1,\cdots ,\omega_m)$. Set $\nu_{\vec{\omega}}=\prod_{i=1}^m \omega_i^{\frac{p}{p_i}}$. We say $\vec{\omega}$ satisfies the $A_{\vec{P}}$ condition if
$$
[\vec{\omega}]_{A_{\vec{P}}}:=\sup_Q \prod_{i=1}^m \left(\frac{1}{|Q|} \int _Q \omega_i^{1-p_i^\prime}\right)^{\frac{p}{p_i^\prime}}  \left(\frac{1}{|Q|} \int_Q \nu_{\vec{\omega}}\right) <\infty,
$$
where $[\vec{\omega}]_{A_{\vec{P}}}$ is called the $A_{\vec{P}}$ constant of $\vec{\omega}$. When $p_i=1$, $\left(\frac{1}{|Q|} \int _Q \omega_i^{1-p_i^\prime}\right)^{\frac{p}{p_i^\prime}}$ is understood as $(\inf_Q \omega_i)^{-1}$.
\end{definition}
 As is well known, the classical $A_p$ and $A_{p,q}$ weights are quite different, so the definition of $A_{(\vec{P},q)}$ weights can not follow the
way to define $A_{\vec{P}}$ weights as in Definition \ref{def 1.2}. The following multiple weights class $A_{(\vec{P},q)}$ was first defined by Chen, Xue in \cite{CX}, and also simultaneously by Moen in \cite{KM}.
\begin{definition} \label{def 1.3} (Class of $A_{(\vec{P},q)}$) (See \cite{KM} or \cite{CX}) Let
$1< p_1,\cdots ,p_m<\infty,\,  \frac{1}{p}=\frac{1}{p_1}+\cdots+\frac{1}{p_m},\,$ and $q>0.$
Suppose that $\vec{\omega}=(\omega_1,\cdots ,\omega_m),\,$ and each $\omega_i$ $(i=1,\cdots ,m)$ is a nonnegative function on $\R^n$.
We say $\vec{\omega}\in A_{(\vec{P},q)}$, if it satisfies
$$
\sup_Q \left(\frac{1}{|Q|} \int _ Q \nu_{\vec{\omega}}^q\right)^{\frac{1}{q}} \prod_{i=1}^m
\left(\frac{1}{|Q|} \int _ Q \omega_i^{-p_i^\prime}\right)^{\frac{1}{p_i^\prime}}<\infty,
$$
where $\nu_{\vec{\omega}} = \prod_{i=1}^m \omega_i$.
If $p_i=1$, $(\frac{1}{|Q|} \int _ Q \omega_i^{-p_i^\prime})^{\frac{1}{p_i^\prime}}$ is understood as $(\inf_Q \omega_i)^{-1}$.
And we denote
$$
[\vec{\omega}]_{A_{(\vec{P},q)}}=\sup_Q \left(\frac{1}{|Q|} \int_Q \nu_{\vec{\omega}}^q\right) \prod_{i=1}^m
\left(\frac{1}{|Q|} \int _ Q \omega_i^{-p_i^\prime}\right)^{\frac{q}{p_i^\prime}}.
$$
\end{definition}
It's easy to see that in the linear case ($m=1$), $A_{(\vec{P},q)}$ will be degenerated to the classical weight class $A_{p,q}$ with $[\vec{\omega}]_{A_{(\vec{P},q)}}=[\omega]_{A_{p,q}}$. Thus $A_{(\vec{P},q)}$ is a natural $m$-linear generalization of the classical $A_{p,q}$ weights which was defined and studied in \cite{MW}.

On one hand, in \cite{LOPTT}, the following multilinear extension of Muckenhoupt $A_p$ theorem for the maximal function was obtained: the inequality
$$
||\mathcal{M}(\vec{f})||_{L^p(\nu_{\vec{\omega}})} \leq C \prod_{i=1}^m ||f_i||_{L^{p_i}(\omega_i)}
$$
holds for every $\vec{f}$ if and only if $\vec{\omega}\in A_{\vec{P}}$, where $\nu_{\vec{\omega}}=\prod_{i=1}^m \omega_i^{\frac{p}{p_i}}$. Therefore, an interesting problem arises naturally, that is: Can we extend Lacy et al's results, inequalities (\ref{1.3}) and (\ref{1.4}) to the multilinear case? For $\alpha=0$, efforts have been made by Dami\'{a}n, Lerner and P\'{e}rez \cite{AKL2}, where they extended Buckley's result \cite{SB} (inequality (\ref{1.3}) for $\alpha=0$) to the multilinear case. More precisely, they proved the following mixed $A_p-A_{\infty}$ estimates, which is sharp in the sense that the exponents can not be replaced by smaller ones.

\vspace{0.2cm}
\noindent\textbf{Theorem A}  {\rm (\cite{AKL2}).}
Let $1<p_i<\infty,\,i=1,\cdots ,m$ and $\frac{1}{p}=\sum_{i=1}^m \frac{1}{p_i}$. Then the inequality

$$
||\mathcal{M}(\vec{f})||_{L^p(\nu_{\vec{\omega}})}\leq C_{n,m,\vec{P}} [\vec{\omega}]_{A_{\vec{P}}}^{1/p} \prod_{i=1}^m [\sigma_i]_{A_\infty}^{1/p_i} \prod_{i=1}^m ||f_i||_{L^{p_i}(\omega_i)}
$$
holds if $\vec{\omega} \in A_{\vec{P}}$, where $\nu_{\vec{\omega}}=\prod_{i=1}^m \omega_i^{\frac{p}{p_i}}$ and $\sigma_i=\omega_i^{1-p_i^\prime},\,i=1,\cdots ,m$.

Recently, Li, Moen and Sun \cite{KKW} obtained the following $A_p$ type estimate which improved the result in \cite{AKL2}.

\vspace{0.2cm}
\noindent\textbf{Theorem B} {\rm (\cite{KKW}).}
Let $1<p_i<\infty ,\, i=1,\cdots ,m$ and $\frac{1}{p}=\sum_{i=1}^m \frac{1}{p_i}$. Denote by $\gamma=\gamma(p_1,\cdots ,p_m)$ the possible best power in the following inequality
$$
||\mathcal{M}(\vec{f})||_{L^p(\nu_{\vec{\omega}})}\leq C_{n,m,\vec{P}} [\vec{\omega}]_{A_{\vec{P}}}^\gamma \prod_{i=1}^m ||f_i||_{L^{p_i}(\omega_i)},
$$
where $\nu_{\vec{\omega}}=\prod_{i=1}^m \omega_i^{\frac{p}{p_i}}$. Then $\frac{m}{mp-1}\leq \gamma \leq \max\{p_1^\prime/p,\cdots ,p^\prime/p\}$.

On the other hand, the multilinear fractional maximal type operator has been studied by Chen, Xue \cite{CX}, and also simultaneously by Moen \cite{KM}. We summarize some weighted norm inequalities for $\mathcal{M}_\alpha$ as follows.

\vspace{0.2cm}
\noindent\textbf{Theorem C} (\cite{KM} or \cite{CX})
Let $0< \alpha<mn,\,1<p_i<\infty,\, i=1,\cdots ,m ,\, \frac{1}{p}=\frac{1}{p_1}+\cdots+\frac{1}{p_m},\, \vec{P}=(p_1,\cdots ,p_m),\, \frac{1}{m}<p<\frac{n}{\alpha} ,\, \frac{1}{q}=\frac{1}{p}-\frac{\alpha}{n}$, and $\nu_{\vec{\omega}} = \prod_{i=1}^m \omega_i$. Then $\vec{\omega}\in A_{(\vec{P},q)}$, if and only if
$$||\mathcal{M}_\alpha(\vec{f})||_{L^q(\nu_{\vec{\omega}}^q)} \leq C \prod\limits_{i=1}^m ||f_i ||_{L^{p_i}(\omega_i^{p_i})}.$$

\vspace{0.2cm}
\noindent\textbf{Theorem D} (\cite{KM} or \cite{CX})
Let $0\le \alpha<mn,\,1\leq p_i<\infty,\, i=1,\cdots ,m ,\, \frac{1}{p}=\frac{1}{p_1}+\cdots+\frac{1}{p_m},\, \vec{P}=(p_1,\cdots ,p_m),\, \frac{1}{m}<p<\frac{n}{\alpha} ,\, \frac{1}{q}=\frac{1}{p}-\frac{\alpha}{n}$, and $\nu_{\vec{\omega}} = \prod_{i=1}^m \omega_i$. Then for $\vec{\omega}\in A_{(\vec{P},q)}$, there is a constant $C>0$ independent of $\vec{f}$ such that $$||\mathcal{M}_\alpha(\vec{f})\nu_{\vec{\omega}}||_{L^{q,\infty}(\R^n)} \leq C \prod\limits_{i=1}^m ||f_i ||_{L^{p_i}(\omega_i^{p_i})}.$$

The first main purpose of this paper is to extend Lacy's result inequality (\ref{1.3}) to the multilinear case for $\alpha>0$. Indeed, we can prove even much more better results with pretty much rough kernels. Before stating our results, we need to introduce the definition of the multilinear fractional maximal type operator with rough homogeneous kernels.
\begin{definition} \label{def 1.4}($\mathcal{M}_{\Omega,\alpha}$ with homogeneous kernels) Assume $\Omega\in L^s((S^{n-1})^m)$ $(s>1)$ is a homogeneous function with degree zero on $\mathbb{R}^n$, i.e. for any $\lambda > 0$ and $y_1,\cdots ,y_m\in \mathbb{R}^n$, $\Omega(\lambda y_1,\cdots ,\lambda y_m)=\Omega(y_1,\cdots ,y_m)$. Then for any $x\in \mathbb{R}^n$ and $0\le \alpha<mn $, define the multilinear fractional maximal type operator with rough homogeneous kernels by
\begin{align}
\mathcal{M}_{\Omega,\alpha}(\vec{f})(x)=\sup_{Q\ni x}  \frac{1}{|Q|^{m-\frac{\alpha}{n}}} \int_{Q^m} |\Omega(x-y_1,\cdots ,x-y_m) |\prod_{i=1}^m |f_i(y_i)| \; d\vec{y},\label{e 1.5}
\end{align}where the supremum in (\ref{e 1.5}) is taken over all cubes $Q$ containing $x$ in $\R^n$ with the sides parallel to the axes. We simply denote $\mathcal{M}_{\Omega,0}$ by $\mathcal{M}_{\Omega}$.
\end{definition}
\remark\label{rem 1.5}
 If we choose $\Omega(y_1,\cdots ,y_m)=1$ in (\ref{e 1.5}), then this operator coincides with the one defined in (\ref{e 1.1}). And if we choose $\Omega(y_1,\cdots ,y_m)=\prod_{i=1}^m \Omega_i(y_i)$, we obtain the operator considered in \cite{CX}.

We obtain a mixed $A_{(\vec{P},q)}-A_\infty$ bound and a $A_{(\vec{P},q)}$ estimate for $\mathcal{M}_{\Omega,\alpha}$ in weighted $L^p$ spaces with different weights as follows:
\begin{theorem}\label{thm 1.1}
Let $1\leq s^\prime<p_i<\infty,\, i=1,\cdots,m$, and $\frac{1}{p}=\sum_{i=1}^m \frac{1}{p_i},\, \vec{P}=(p_1,\cdots ,p_m), \, 0<\alpha<mn, \,\frac{1}{m}<p<\frac{n}{\alpha} ,\, \frac{1}{q}=\frac{1}{p}-\frac{\alpha}{n}$, and $\nu_{\vec{\omega}} = \prod_{i=1}^m \omega_i$. Define $(\omega_1^{s^\prime},\cdots ,\omega_m^{s^\prime})$ and $(p_1/s^\prime,\cdots ,p_m/s^\prime)$ as $\vec{\omega}^{s^\prime}$ and $\frac{\vec{P}}{s^\prime}$, respectively. Assume
$\vec{\omega}^{s^\prime}\in A_{(\frac{\vec{P}}{s^\prime},\frac{q}{s^\prime})}$,
then there exists a constant $C=C(m,n,\vec{P},q,s)$, such that
\begin{align}\label{e 1.6}
||\mathcal{M}_{\Omega,\alpha}(\vec{f})||_{L^q(\nu_{\vec{\omega}}^q)} \leq C ||\Omega||_{L^s((S^{n-1})^m)} [\vec{\omega}^{s^\prime}]_{A_{(\frac{\vec{P}}{s^\prime},\frac{q}{s^\prime})}}^{\frac{1}{q}} \prod _{i=1}^m [\omega_i^{-(\frac{p_i}{s^\prime})^\prime}]_{A_\infty}^{\frac{1}{p_i}} \prod_{i=1}^m
||f_i||_{L^{p_i}(\omega_i^{p_i})}.
\end{align}

\end{theorem}

\begin{theorem}\label{thm 1.2}
Let $1\leq s^\prime<p_i<\infty,\, i=1,\cdots,m$, and $\frac{1}{p}=\sum_{i=1}^m \frac{1}{p_i},\, \vec{P}=(p_1,\cdots ,p_m), \, 0<\alpha<mn, \,\frac{1}{m}<p<\frac{n}{\alpha} ,\, \frac{1}{q}=\frac{1}{p}-\frac{\alpha}{n}$, and $\nu_{\vec{\omega}} = \prod_{i=1}^m \omega_i$. Define $(\omega_1^{s^\prime},\cdots ,\omega_m^{s^\prime})$ and $(p_1/s^\prime,\cdots ,p_m/s^\prime)$ as $\vec{\omega}^{s^\prime}$ and $\frac{\vec{P}}{s^\prime}$, respectively. Define by $\gamma=\gamma(p_1,\cdots ,p_m)$ the power in
\begin{align}\label{e 1.7}
||\mathcal{M}_{\Omega,\alpha}(\vec{f})||_{L^q(\nu_{\vec{\omega}}^q)} \leq C_{n,m,\vec{P},q,s} ||\Omega||_{L^s((S^{n-1})^m)} [\vec{\omega}^{s^\prime}]_{A_{(\frac{\vec{P}}{s^\prime},\frac{q}{s^\prime})}}^\gamma ||f_i||_{L^{p_i}(\omega_i^{p_i})}.
\end{align}
Then for all $1\leq s^\prime<p_1,\cdots ,p_m<\infty,\, $ we have
$$\frac{mp}{q(mp-s^\prime)}(1-\frac{s^\prime\alpha}{mn})\leq \gamma \leq \frac{1}{q}(1-\frac{s^\prime\alpha}{mn})\max\{(\frac{p_1}{s^\prime})^\prime,\cdots ,(\frac{p_m}{s^\prime})^\prime\}.$$
Especially, if $1\leq s^\prime < p_1=\cdots=p_m=s^\prime r<\infty$, we have $\gamma=r^\prime(1-\frac{s^\prime\alpha}{mn})/q$.
\end{theorem}
Take $s^\prime=1$ and $\Omega(y_1,\cdots ,y_m)=1$ in the above theorem, we obtain the following corollaries for $\mathcal{M}_\alpha$. As was shown in (\cite{AKL2}), the first corollary is sharp when $\alpha=0$.
\begin{corollary} \label{cor 1.3}

Let $1<p_i<\infty,\, i=1,\cdots ,m$, and $\frac{1}{p}=\frac{1}{p_1}+\cdots+\frac{1}{p_m},\, \vec{P}=(p_1,\cdots ,p_m),\,  0<\alpha<mn,\, \frac{1}{m}<p<\frac{n}{\alpha}$ and $\frac{1}{q}=\frac{1}{p}-\frac{\alpha}{n}$. Then the inequality
\begin{equation}\label{e 1.8}
||\mathcal{M}_\alpha(\vec{f})||_{L^q(\nu_{\vec{\omega}}^q)} \leq C_{n,m,\vec{P},q} [\vec{\omega}]^{\frac{1}{q}}_{A_{(\vec{P},q)}} \prod_{i=1}^m [\sigma_i]_{A_{\infty}}^{\frac{1}{p_i}} \prod_{i=1}^m ||f_i ||_{L^{p_i}(\omega_i^{p_i})},
\end{equation}
holds if $\vec{\omega}\in A_{(\vec{P},q)}$, where $\nu_{\vec{\omega}} = \prod_{i=1}^m \omega_i, \, \sigma_i=\omega_i^{-p_i^\prime}, \; (i=1,\cdots ,m)$.

\end{corollary}

\begin{corollary} \label{cor 1.4}

Let $1<p_i<\infty,\,i=1,\cdots ,m$, and $\frac{1}{p}=\frac{1}{p_1}+\cdots+\frac{1}{p_m},\, \vec{P}=(p_1,\cdots ,p_m), \, 0<\alpha<mn, \,\frac{1}{m}<p<\frac{n}{\alpha}$ and $\frac{1}{q}=\frac{1}{p}-\frac{\alpha}{n}$. Define by $\gamma=\gamma(p_1,\cdots ,p_m)$ the power in
\begin{align}\label{e 1.9}
||\mathcal{M}_\alpha(\vec{f})||_{L^q(\nu_{\vec{\omega}}^q)} \leq C_{n,m,\vec{P},q} [\vec{\omega}]_{A_{(\vec{P},q)}}^\gamma \prod_{i=1}^m ||f_i ||_{L^{p_i}(\omega_i^{p_i})},
\end{align}
where $\nu_{\vec{\omega}}$ is the same as in Theorem \ref{thm 1.9}.
Then for all $1<p_1,\cdots ,p_m<\infty, \,$ we have $$\frac{mp}{q(mp-1)}(1-\frac{\alpha}{mn}) \leq \gamma \leq \frac{1}{q}(1-\frac{\alpha}{mn})\max \{p_1^\prime,\cdots ,p_m^\prime\}.$$ Especially, if $1<p_1=\cdots=p_m=r<\infty$, we have $\gamma=\frac{r\prime}{q}(1-\frac{\alpha}{mn})$. Furthermore, the power is sharp.
\remark Note that if $m=1$ in Corollary \ref{cor 1.4}, then the sharp power $\gamma=\frac {p'}q(1-\frac \alpha n)$. This coincides with (\ref{1.3}). The above results are new for $m\ge 2.$

\end{corollary}

We also obtain a weak type estimate as follows,
\begin{theorem}\label{thm 1.5}
Let $0<\alpha<mn,\,1\leq p_i<\infty,\, i=1,\cdots ,m ,\, \frac{1}{p}=\frac{1}{p_1}+\cdots+\frac{1}{p_m},\, \vec{P}=(p_1,\cdots ,p_m),\, \frac{1}{m}<p<\frac{n}{\alpha} ,\, \frac{1}{q}=\frac{1}{p}-\frac{\alpha}{n}$, and $\nu_{\vec{\omega}} = \prod_{i=1}^m \omega_i$. Then for $\vec{\omega}\in A_{(\vec{P},q)}$, there is a constant $C>0$ independent of $\vec{f}$ such that
\begin{align}\label{e 1.15}||\mathcal{M}_\alpha(\vec{f})\nu_{\vec{\omega}}||_{L^{q,\infty}(\R^n)} \leq C [\vec{\omega}]_{A_{(\vec{P},q)}}^{\frac1 q}\prod\limits_{i=1}^m ||f_i ||_{L^{p_i}(\omega_i^{p_i})}.
\end{align}
\end{theorem}

\remark
when $m=1$, $\alpha=0$ this can be found in Garc¨ªa-Cuerva and Rubio de Francia's book. If $m=1$, $\alpha>0$ it was given in \cite{LMPT} by Lacy et al.

Corollary \ref{cor 1.4} can be formulated in a more general setting, before stating our results, we begin with one more definition.
\begin{definition} \label{def 1.5}
Let $X$ be a Banach function space, $X^\prime$ is the associate space to $X$. Given a cube $Q$, define the $X$-average of $f$ over $Q$ and the maximal operator $M_X$ by
$$
||f||_{X,Q} =||\tau_{l(Q)} (f\chi_Q)||_X,\quad M_Xf(x)= \sup_{Q\ni x} ||f||_{X,Q},
$$
where $l(Q)$ denotes the side length of $Q$ and $\tau_\delta f(x)=f(\delta x), \, \delta>0,\, x\in \R^n$.
\end{definition}
\begin{theorem} \label{thm 1.6}
Let $0<\alpha<mn,\,1<p_i<\infty,\; i=1,\cdots ,m$, and $\frac{1}{p}=\frac{1}{p_1}+\cdots+\frac{1}{p_m}, \frac{1}{m}<p<\frac{n}{\alpha}$ and $\frac{1}{q}=\frac{1}{p}-\frac{\alpha}{n}$. Let $X_i$ be a Banach function space such that $M_{X_i^\prime}$ is bounded on $L^{p_i}(\R^n)$. Let $u$ and $v_1,\cdots ,v_m$ be the weights satisfying
$$
K=\sup_Q \left(\frac{u(Q)}{|Q|}\right)^{\frac{1}{q}} \prod_{i=1}^m ||v_i^{-1}||_{X_i,Q} < \infty.
$$
Then
$$
||\mathcal{M}_\alpha(\vec{f})||_{L^q(u)} \leq CK \prod _{i=1}^m ||M_{X_i^\prime}||_{L^{p_i}(\R^n)} ||f_i ||_{L^{p_i}(v_i^{p_i})}.
$$
\end{theorem}
As applications of the above results, we consider the multilinear fractional integral operators with rough kernels defined by
\begin{equation}\mathcal{I}_{\Omega,\alpha}(\vec{f})(x) = \int_{{(\mathbb{R}^n)}^m}\frac{\Omega(y_1,\cdots, y_m)}{{|(y_1, \cdots, y_m)|}^{mn - \alpha}}
\prod_{i = 1}^{m}f_i(x - y_i)\,d\vec{y},\end{equation}
we have some almost sharp estimates for the fractional integral operators as follows:
\begin{corollary}[Weighted strong bounds of $\mathcal{I}_{\Omega, \alpha}$]\label{cor 1.7}
Let $m\ge 2,\, 0 < \alpha < mn$, $1 \leqslant s' < p_1, \cdots, p_m <
\infty$, $\frac{1}{p} = \frac{1}{p_1} + \cdots + \frac{1}{p_m}$,
and $\frac{1}{q} =
\frac{1}{p} - \frac{\alpha}{n}$. Denote $\vec{\omega}^{s'} =
(\omega_1^{s'}, \cdots, \omega_m^{s'})$ and $\frac{\vec{P}}{s'} =
(\frac{p_1}{s'}, \cdots, \frac{p_m}{s'})$. Assume
$\vec{\omega}^{s'}\in A_{(\frac{\vec{P}}{s'}, \frac{q}{s'})} \cap
A_{(\frac{\vec{P}}{s'}, \frac{q_{\epsilon}}{s'})} \cap
A_{(\frac{\vec{P}}{s'}, \frac{q_{- \epsilon}}{s'})}$, where $0<\epsilon <\alpha$,
$\frac{1}{q_{\epsilon}} = \frac{1}{p} - \frac{\alpha + \epsilon}{n}$
and $\frac{1}{q_{- \epsilon}} = \frac{1}{p} - \frac{\alpha -
\epsilon}{n}$. Then there
is a constant $C
> 0$ independent of $\vec{f}$ such that
\begin{equation}\label{Weighted Strong kernel}
{\big\|\mathcal{I}_{\Omega,
\alpha}(\vec{f})\big\|}_{L^q({\nu_{\vec{\omega}}}^{q})} \leqslant C_{||\Omega||_{L^s} } \bigg\{[\vec{\omega}^{s^\prime}]_{A_{(\frac{\vec{P}}{s^\prime},\frac{q_\epsilon}{s^\prime})}}^{\frac{1}{2q_\epsilon}}[\vec{\omega}^{s^\prime}]_{A_{(\frac{\vec{P}}{s^\prime},\frac{q_{-\epsilon}}{s^\prime})}}^{\frac{1}{2q_{-\epsilon}}}\prod _{i=1}^m [\omega_i^{-(\frac{p_i}{s^\prime})^\prime}]_{A_\infty}^{\frac{1}{p_i}}\bigg\}
\prod_{i =
1}^m{\big\|f_i\big\|}_{L^{p_i}(\omega_i^{p_i})}.\end{equation}
Moreover, if $s'=1$ and $\Omega\equiv 1$, then we have \begin{equation}\label{Weighted Strong kernel}
{\big\|\mathcal {I}_
\alpha(\vec{f})\big\|}_{L^q({\nu_{\vec{\omega}}}^{q})} \leqslant C \bigg\{[\vec{\omega}]_{A_{(\vec{P},q_{\epsilon})}}^{\gamma_1}[\vec{\omega}]_{A_{(\vec{P},q_{-\epsilon})}}^{\gamma_2} \bigg\}^{1/2}\prod _{i=1}^m
{\big\|f_i\big\|}_{L^{p_i}(\omega_i^{p_i})}\end{equation} where
$\frac{mp}{q_{\epsilon}(mp-1)}(1-\frac{\alpha+\epsilon}{mn}) \leq \gamma_1 \leq \frac{1}{q_{\epsilon}}(1-\frac{\alpha+\epsilon}{mn})\max \{p_1^\prime,\cdots ,p_m^\prime\}$ and $\frac{mp}{q_{-\epsilon}(mp-1)}(1-\frac{\alpha-\epsilon}{mn}) \leq \gamma_2 \leq \frac{1}{q_{-\epsilon}}(1-\frac{\alpha-\epsilon}{mn})\max \{p_1^\prime,\cdots ,p_m^\prime\}$, Especially, if $1<p_1=\cdots=p_m=r<\infty$, we have $\gamma_1=\frac{r\prime}{q_{\epsilon}}(1-\frac{\alpha+\epsilon}{mn})$ and $\gamma_2=\frac{r\prime}{q_{-\epsilon}}(1-\frac{\alpha-\epsilon}{mn}).$

\end{corollary}

\begin{remark}The almost sharp estimates come from the fact that, if $\epsilon =0$, $m=1$ and $\frac {p'}q>1$, we can obtain $\gamma =\frac {p'}{q}(1-\frac \alpha n)$. In this case, the estimate in Corollary \ref{cor 1.7} is sharp. In the linear case, it is well known that
\cite[pp. 152]{Singular Integral and Related Topics}
$\omega^{s'}\in A_{(\frac{p}{s'}, \frac{q}{s'})}$ concludes
$\omega^{s'}\in A_{(\frac{p}{s'}, \frac{q_{\epsilon}}{s'})} \cap
A_{(\frac{p}{s'}, \frac{q_{- \epsilon}}{s'})}$ by the monotonicity
of class $A_p$. However, unlike $A_{p, q}$, $A_{(\vec{P}, q)}$ doesn't have
such good property because of $A_{\vec{P}}$ is not monotone. Therefore, we have to assume that
$\vec{\omega}^{s'} \in A_{(\frac{\vec{P}}{s'},
\frac{q_{\epsilon}}{s'})} \cap A_{(\frac{\vec{P}}{s'}, \frac{q_{-
\epsilon}}{s'})}$.\end{remark}

One may ask if the above results still hold for $\mathcal{M}_\Omega$ or not, we summarize some results as follows:

\begin{theorem}\label{thm 1.8}(Mixed $A_p-A_\infty$ estimates)
Let $1\leq s^\prime<p_i<\infty,\, i=1,\cdots ,m$, and  $\frac{1}{p}=\sum_{i=1}^m \frac{1}{p_i}$. Then for $\vec{\omega} \in A_{\vec{P}/s^\prime}$, we have
\begin{align}\label{e 1.13}
||\mathcal{M}_\Omega(\vec{f})||_{L^p(\nu_{\vec{\omega}})}\leq C_{n,m,\vec{P},s} ||\Omega||_{L^s((S^{n-1})^m)}  [\vec{\omega}]_{A_{\vec{P}/s^\prime}}^{1/p} \prod_{i=1}^m [\omega_i^{1-(\frac{p_i}{s^\prime})^\prime}]_{A_\infty}^{1/p_i}\prod_{i=1}^m ||f_i||_{L^{p_i}(\omega_i)},
\end{align}
where $\nu_{\vec{\omega}} = \prod_{i=1}^m \omega_i^{\frac{p}{p_i}},\, i=1,\cdots ,m$.

\end{theorem}
\begin{theorem}\label{thm 1.9}($A_p$ estimates)
Let $1\leq s^\prime<p_i<\infty,\, i=1,\cdots ,m$, and  $\frac{1}{p}=\sum_{i=1}^m \frac{1}{p_i}$. Denote by $\gamma=\gamma(p_1,\cdots ,p_m)$ the power in the following inequality

\begin{align}\label{e 1.14}
||\mathcal{M}_\Omega(\vec{f})||_{L^p(\nu_{\vec{\omega}})}\leq C_{n,m,\vec{P},s} ||\Omega||_{L^s((S^{n-1})^m)}   [\vec{\omega}]_{A_{\vec{P}/s^\prime}}^\gamma \prod_{i=1}^m ||f_i||_{L^{p_i}(\omega_i)},
\end{align}
where $\nu_{\vec{\omega}} = \prod_{i=1}^m \omega_i^{\frac{p}{p_i}}$. Then
$\frac{m}{mp-s^\prime} \leq \gamma \leq \max\{\frac{(p_1/s^\prime)^\prime}{p},\cdots ,\frac{(p_m/s^\prime)^\prime}{p}\}$.

\end{theorem}
\remark\label{rem 1.2}
Theorem \ref{thm 1.8} and Theorem \ref{thm 1.9} are the generalization of Theorem A and Theorem B. In fact, we only need to choose $s^\prime=1$ and $\Omega(y_1,\cdots ,y_m)=1$ in (\ref{e 1.13}) and (\ref{e 1.14}).

In this paper, we will prove Theorem \ref{thm 1.1} and Corollary \ref{cor 1.3} in Section 2. In Section 3, we give the proof of Theorem \ref{thm 1.2}, Corollary \ref{cor 1.4} and Theorem \ref{thm 1.5}. We will prove a two weights version of (\ref{e 1.8}) for a general Banach Space Theorem \ref{thm 1.6} in Section 4. In the last part, we show the proof of the Corollary \ref{cor 1.7}, Theorem \ref{thm 1.8} and Theorem \ref{thm 1.9}.

\section{The proof of Theorem \ref{thm 1.1} and Corollary \ref{1.3}}
\setcounter{equation}{0}

To begin with, we prepare two definitions. Recall that the standard dyadic grid in $\R^n$ consists of the cubes
$$
2^{-k}([0,1)^n+j),\, k\in \mathbb{Z},\, j\in \mathbb{Z}^n
$$
\begin{definition} \label{def 2.1}
We define a general dyadic grid $\mathcal{D}$, if it is a collection of cubes with the following properties:\\
$(i)$for any $Q\in \mathcal{D}$, its sidelength $l(Q)$ is of the form $2^k, k\in \mathbb{Z}$;\\
$(ii)Q\cap R\in \{Q,R,\emptyset\}$ for any $Q,R\in \mathcal{D}$;\\
$(iii)$the cubes of a fixed sidelength $2^k$ form a partition of $\R^n$.
\end{definition}
We define $\mathcal{D}_\beta= \{2^{-k}([0,1)^n+j+(-1)^k \beta),\, k\in \mathbb{Z},\, j\in \mathbb{Z}^n\},  \, \beta\in \{0,\frac{1}{3}\}^n$.\\
\begin{definition} \label{def 2.2}
We say that $\{Q_j^k\}$ is a sparse family of cubes, if\\
$(i)$the cubes $Q_j^k$ are disjoint in $j$, with $k$ fixed;\\
$(ii)$if $\Omega _k= \bigcup _j Q_j^k,\,$ then $\Omega _{k+1}\subseteq \Omega _k$;\\
$(iii)|\Omega _{k+1}\cap Q_j^k|\leq \frac{1}{2} |Q_j^k|$.
\end{definition}
\remark \label{rem 2.3}
With each sparse family $\{Q_j^k\}$, we associate the sets $E_j^k=Q_j^k\backslash \Omega_{k+1}$. Observe that the sets $E_j^k$ are pairwise disjoint, $\cup_{j,k} E_j^k= \R^n$and $|Q_j^k|\leq 2|E_j^k|$.\\

We will use the following proposition in \cite{THCP}, and the proof can be found in \cite{AKL1}.

\begin{proposition} \label{prop 2.1}
There are $2^n$ dyadic grids $\mathcal{D}_\beta$ such that for any cube $Q\subset \R^n$, there exists a cube $Q_\beta \in \mathcal{D}_\beta$ such that $Q\subseteq Q_\beta$ and $l(Q_\beta)\leq 6 l(Q)$.
\end{proposition}
We need the following reverse H\"{o}lder property of $A_\infty$ weights which was proved in \rm{(\cite{THCP})}:
\begin{lemma} \rm{(\cite{THCP})}. Let $\omega\in A_\infty$, then the following inequality holds
$$
\left(\frac{1}{|Q|} \int _Q \omega^{r(\omega)}\right)^{\frac{1}{r(\omega)}} \leq 2 \frac{1}{|Q|} \int_Q \omega,
$$
where $r(\omega)=1+\frac{1}{c_n [\omega]_{A_\infty}}$, $c_n=2^{11+n}$ and it is easy to see $r(\omega)^\prime \approx [\omega]_{A_\infty}$.\end{lemma}

\noindent
{\textbf{\it The proof of Corollary \ref{1.3}.}}

First, by Proposition \ref{prop 2.1}, for any $Q\subset \R^n$,

\begin{align*}
\prod_{i=1}^m &\frac{1}{|Q|^{1-\frac{\alpha}{mn}}} \int _Q |f_i(y_i)|\; dy_i\nonumber\le \prod_{i=1}^m \frac{6^{n(1-\frac{\alpha}{mn})}}{|Q_\beta|^{1-\frac{\alpha}{mn}}} \int _{Q_\beta} |f_i(y_i)|\; dy_i\nonumber\\
\leq& 6^{mn(1-\frac{\alpha}{mn})} \sum_{\beta\in \{0,\frac{1}{3}\}^n} \prod _{i=1}^m \frac{1}{|Q_\beta|^{1-\frac{\alpha}{mn}}} \int _{Q_\beta} |f_i(y_i)|\; dy_i\nonumber
\leq 6^{mn-\alpha} \sum_{\beta\in \{0,\frac{1}{3}\}^n} \mathcal{M}_\alpha^{\mathcal{D}_\beta}(\vec{f})(x)
\end{align*}

So we obtain
\begin{align}
\mathcal{M}_\alpha(\vec{f})(x) \leq 6^{mn-\alpha} \sum_{\beta\in \{0,\frac{1}{3}\}^n} \mathcal{M}_\alpha^{\mathcal{D}_\beta}(\vec{f})(x).\label{e 2.1}
\end{align}
By (\ref{e 2.1}), it suffices to prove Corollary \ref{1.3} for the dyadic multilinear fractional maximal type operator $\mathcal{M}_\alpha^{\mathcal{D}_\beta}$. Since the proof is independent of the particular dyadic grid, without loss of generality, we consider $\mathcal{M}_\alpha^d$ taken with respect to the standard dyadic grid. We will use exactly the argument as in the Calder\'{o}n-Zygmund decomposition.
For a real number $a$, which will be specified below and for $k\in \mathbb{Z}$, consider the sets $\Omega_k=\{x\in \R^n :\mathcal{M}_\alpha^d(\vec{f})(x)>a^k\}$.
Then we have that $\Omega_k=\cup_j Q_j^k$, where the cubes $Q_j^k$ are pairwise disjoint with $k$ fixed, and
$$
a^k<\prod_{i=1}^m \frac{1}{|Q_j^k|^{1-\frac{\alpha}{mn}}} \int _{Q_j^k} |f_i(y_i)|\; dy_i \leq 2^{mn} a^k.
$$
From this and according to H\"{o}lder's inequality,
\begin{align*}
|Q_j^k\cap \Omega_{k+1}|=& \sum_{Q_l^{k+1}\subset Q_j^k} |Q_l^{k+1}|\le a^{-\frac{k+1}{m}} \sum_{Q_l^{k+1}\subset Q_j^k} \prod_{i=1}^m \left(|Q_l^{k+1}|^{\frac{\alpha}{mn}}\int _{Q_l^{k+1}} |f_i(y_i)| \; dy_i \right)^{\frac{1}{m}}\\
\leq &a^{-\frac{k+1}{m}} \prod_{i=1}^m \left(|Q_j^k|^{\frac{\alpha}{mn}}\int _{Q_j^k} |f_i(y_i)| \; dy_i \right)^{\frac{1}{m}}
\leq 2^n a^{-\frac{1}{m}} |Q_j^k|.
\end{align*}
Hence, taking $a = 2^{m(n+1)}$, we obtain that the family $\{Q_j^k\}$ is sparse. Then we can see
\begin{align*}
\int _{\mathbb{R}^n} |\mathcal{M}_\alpha^d&(\vec{f})(x)\nu_{\vec{\omega}}(x)|^q \; dx= \sum_k \int_{\Omega_k\backslash\Omega_{k+1}} |\mathcal{M}_\alpha^d(\vec{f})(x)\nu_{\vec{\omega}}(x)|^q \; dx\\
\leq& a^q \sum_{k,j} \left(\prod_{i=1}^m \frac{1}{|Q_j^k|^{1-\frac{\alpha}{mn}}} \int _{Q_j^k} |f_i| \; dy_i\right)^q \int_{Q_j^k} \nu_{\vec{\omega}}(x)^q \; dx\\
\leq&  a^q \sum_{k,j} \left(\prod_{i=1}^m \frac{1}{|Q_j^k|} \int _{Q_j^k} |f_i \omega_i|^{\theta_i} \; dy_i\right)^{\frac{q}{\theta_i}} \left(\prod_{i=1}^m \frac{1}{|Q_j^k|} \int _{Q_j^k}  \omega_i^{-\theta_i^\prime} \; dy_i\right)^{\frac{q}{\theta_i^\prime}} \\&\quad\times\frac{|Q_j^k|^{1+\frac{\alpha q}{n}}}{|Q_j^k|} \int_{Q_j^k} \nu_{\vec{\omega}}(x)^q \; dx\\
\leq& C a^q \sum_{k,j} \prod_{i=1}^m \left(\frac{1}{|Q_j^k|} \int _{Q_j^k} |f_i \omega_i|^{\theta_i} \; dy_i\right)^{\frac{q}{\theta_i}} \prod_{i=1}^m\left( \frac{1}{|Q_j^k|} \int _{Q_j^k} \omega_i^{-p_i^\prime} \; dy_i\right)^{\frac{q}{p_i^\prime}}\\&\quad \times \frac{\nu_{\vec{\omega}}^q(Q_j^k)}{|Q_j^k|} |Q_j^k|^{\frac{q}{p}}\\ \leq& C [\vec{\omega}]_{A_{(\vec{P},q)}} \sum_{k,j} \left(\prod_{i=1}^m (\frac{1}{|Q_j^k|} \int _{Q_j^k} |f_i \omega_i|^{\theta_i} \; dy_i)^{\frac{p}{\theta_i}} |E_j^k|\right)^{\frac{q}{p}}\triangleq C [\vec{\omega}]_{A_{(\vec{P},q)}} U(f)
 \end{align*}
Now, we are in position to estimate $U(f).$
\begin{align*}
U(f)
\leq& C \left(\sum_{j,k} \int_{E_j^k} \prod_{i=1}^m M((f_i \omega_i)^{\theta_i})(x)^{\frac{p}{\theta_i}} \; dx \right)^{\frac{q}{p}}
\leq C  \left(\int_{\R^n} \prod_{i=1}^m M((f_i \omega_i)^{\theta_i})(x)^{\frac{p}{\theta_i}} \; dx \right)^{\frac{q}{p}}\\
\leq& C\prod _{i=1}^m \left(\int_{\R^n}  M((f_i \omega_i)^{\theta_i})(x)^\frac{p_i}{\theta_i} \; dx \right)^{\frac{q}{p_i}}
\leq C  \prod _{i=1}^m \left(\left(\frac{p_i}{\theta_i}\right)^\prime\right)^{\frac{q}{p_i}} \left(\int_{\R^n} |f_i \omega_i|^{p_i} \; dx\right)^{\frac{q}{p_i}}\\
\leq &C \prod _{i=1}^m (r_i^\prime p_i^\prime)^{\frac{q}{p_i}} ||f_i ||_{L^{p_i}(\omega_i^{p_i})}^q
\leq C \prod _{i=1}^m [\sigma_i]_{A_\infty}^{\frac{q}{p_i}} ||f_i ||_{L^{p_i}(\omega_i^{p_i})}^q,
\end{align*}
where $\theta_i=(r_i p_i^\prime)^\prime$ and $r_i$ is the exponent in the sharp reverse H\"{o}lder's inequality for the weight $\sigma_i=\omega_i^{-p_i^\prime}$ which are in $A_\infty$ for $i=1,\cdots ,m$.
Thus, we have
$$
||\mathcal{M}^d_\alpha(\vec{f})||_{L^q(\nu_{\vec{\omega}}^q)} \leq C_{n,m,\vec{P},q} [\vec{\omega}]^{\frac{1}{q}}_{A_{(\vec{P},q)}} \prod_{i=1}^m [\sigma_i]_{A_{\infty}}^{\frac{1}{p_i}} \prod_{i=1}^m ||f_i ||_{L^{p_i}(\omega_i^{p_i})}.
$$
The proof of Corollary \ref{1.3} is finished.\\

\noindent
\text{\it {Proof of Theorem \ref{thm 1.1}.} }

We claim the following inequality.
\begin{align}\label{e 2.2}
\mathcal{M}_{\Omega,\alpha}(\vec{f})(x) \leq C ||\Omega||_{L^s((S^{n-1})^m)}  (\mathcal{M}_{\alpha s^\prime}(\vec{f}^{s^\prime})(x))^{\frac{1}{s^\prime}}.
\end{align}

\noindent
In fact, suppose $B$ with radius $r$ is the circumscribed ball of the cube $Q$,

\begin{align*}
&\mathcal{M}_{\Omega,\alpha} (\vec{f})(x) \\
\leq &\sup_{Q\ni 0} {|Q|^{\frac{\alpha}{n}}} \left( \frac{1}{|Q|^m} \int_{Q^m} \prod_{i=1}^m|f_i(x-y_i)|^{s^\prime} \; d\vec{y}  \right)^{\frac1{s^\prime}}  \left( \frac{1}{|Q|^m} \int_{B^m} |\Omega(y_1,\cdots ,y_m)|^s\; d\vec{y} \right)^{\frac1s}\\
\leq& C ||\Omega||_{L^s((S^{n-1})^m)} \sup_{Q\ni 0} \prod_{i=1}^m  \left(\frac{1}{|Q|^{1-\frac{\alpha s^\prime}{mn}}} \int_Q |f_i(x-y_i)|^{s^\prime} \; dy_i \right)^{1/s^\prime}\\
=& C ||\Omega||_{L^s((S^{n-1})^m)}  \mathcal{M}_{\alpha s^\prime}(\vec{f}^{s^\prime})(x)^{\frac{1}{s^\prime}}.
\end{align*}

Then from Corollary \ref{cor 1.3}, we have
\begin{align*}
&||\mathcal{M}_{\Omega,\alpha}(\vec{f})||_{L^q(\nu_{\vec{\omega}}^q)}
\leq C  ||\Omega||_{L^s((S^{n-1})^m)}  ||\mathcal{M}_{\alpha s^\prime}(\vec{f}^{s^\prime})||_{L^{q/s^\prime}(\nu_{\vec{\omega}^{s^\prime}}^{q/s^\prime})}^{\frac{1}{s^\prime}}\\
\leq& C ||\Omega||_{L^s((S^{n-1})^m)}  [\vec{\omega}^{s^\prime}]_{A_{(\frac{\vec{P}}{s^\prime},\frac{q}{s^\prime})}}^{1/q} \prod _{i=1}^m [\omega_i^{-(\frac{p_i}{s^\prime})^\prime}]_{A_\infty}^{\frac{s^\prime}{p_i} \frac{1}{s^\prime}} \prod_{i=1}^m
||f_i||_{L^{p_i}(\omega_i^{p_i})}\\
\leq& C ||\Omega||_{L^s((S^{n-1})^m)} [\vec{\omega}^{s^\prime}]_{A_{(\frac{\vec{P}}{s^\prime},\frac{q}{s^\prime})}}^{1/q} \prod _{i=1}^m [\omega_i^{-(\frac{p_i}{s^\prime})^\prime}]_{A_\infty}^{\frac{1}{p_i}} \prod_{i=1}^m
||f_i||_{L^{p_i}(\omega_i^{p_i})}.
\end{align*}
where $\nu_{\vec{\omega}^{s^\prime}}^{q/s^\prime}=\prod_{i=1}^m \omega_i^q$.
We thus complete the proof of Theorem \ref{thm 1.1}.

\section{The proof of Theorem \ref{thm 1.2} and Corollary \ref{1.4} and Theorem \ref{thm 1.5}}

\noindent
\text{\it Proof of Corollary \ref{1.4}.}

First, we will show the low bound of $\gamma$.

Assume that
$n=1,\, 0<\varepsilon<1, \,\omega_i(x)=|x|^{(1-\varepsilon)(1-\frac{1}{p_i})},\, f_i(x)= x^{-1+\varepsilon} \chi_{(0,1)}(x), \, i=1,\cdots ,m.$ Then \\
$$
\nu_{\vec{\omega}}(x)=|x|^{\left(1-\varepsilon\right)(m-\frac{1}{p})}, \, [\vec{\omega}]_{A_{(\vec{P},q)}}\approx \left(\frac{1}{\varepsilon}\right)^{q(m-\frac{1}{p})},
$$

$$
\prod_{i=1}^m  ||f_i ||_{L^{p_i}(\omega_i^{p_i})}=\left(\frac{1}{\varepsilon}\right)^{\frac{1}{p}},\, ||\mathcal{M}_\alpha(\vec{f})||_{L^q(\nu_{\vec{\omega}}^q)} \geq \left(\frac{1}{\varepsilon}\right)^{m+\frac{1}{q}}.
$$

So, by (\ref{e 1.9}), we have

$$
m+\frac{1}{q}\leq \gamma q(m-\frac{1}{p})+\frac{1}{p}.
$$
That is, $\gamma\geq \frac{mp}{q(mp-1)}(1-\frac{\alpha}{mn})$.\\

Now, we will show the upper bound for $\gamma$. By the proof of the Corollary \ref{1.3}, we can replaced $\mathcal{M}_\alpha(\vec{f})$ by $\mathcal{M}_\alpha^d(\vec{f})$ for standard dyadic grid in (\ref{e 1.9}). In fact, we only need to show

$$
||\mathcal{M}_\alpha^d(f_1 \sigma_1,\cdots ,f_m \sigma_m)||_{L^q(\nu_{\vec{\omega}}^q)} \leq C_{m,n,\vec{P},q} [\vec{\omega}]_{A_{(\vec{P},q)}}^{\frac{1}{q}(1-\frac{\alpha}{mn})\max\{p_1^\prime,\cdots ,p_m^\prime\}} \prod_{i=1}^m ||f_i||_{L^{p_i}(\sigma_i)}
$$
where $\sigma_i = \omega_i^{-p_i^\prime},\,i=1,\cdots ,m$.

Let $a=2^{m(n+1)}$ and $\Omega_k=\{x\in \R^n :\mathcal{M}_\alpha^d(f_1 \sigma_1,\cdots ,f_m \sigma_m)(x)>a^k\}$. Then we have $\Omega_k=\cup_j Q_j^k$ and $S=\{Q_j^k\}$ is sparse, satisfying
$$
a^k<\prod_{i=1}^m \frac{1}{|Q|^{1-\frac{\alpha}{mn}}} \int _Q |f_i| \sigma_i\; dy_i \leq 2^{mn} a^k.
$$
Without loss of the generality, we assume $p_1=\min\{p_1,\cdots ,p_m\}$. Then

\begin{align*}
&\int_{\R^n} \mathcal{M}_\alpha^d(f_1 \sigma_1,\cdots ,f_m \sigma_m)(x)^q \nu_{\vec{\omega}}(x)^q \; dx\\
=& \sum_k \int_{\Omega_k\backslash\Omega_{k+1}} \mathcal{M}_\alpha^d(f_1 \sigma_1,\cdots ,f_m \sigma_m)(x)^q\nu_{\vec{\omega}}(x)^q \; dx\\
\leq& a^q \sum_{k,j} \left(\prod_{i=1}^m \frac{1}{|Q_j^k|^{1-\frac{\alpha}{mn}}} \int _{Q_j^k} |f_i|\sigma_i \; dy_i\right)^q \int_{Q_j^k} \nu_{\vec{\omega}}(x)^q \; dx\\
\leq& a^q \sum_{k,j} \left(\frac{\nu_{\vec{\omega}}^q(Q_j^k)^{p_1^\prime}\prod_{i=1}^m \sigma_i(Q_j^k)^{\frac{qp_1^\prime}{p_i^\prime}}}{|Q_j^k|^{[1+(m-\frac{1}{p})q]p_1^\prime}}\right)^{(1-\frac{\alpha}{mn})}  \cdot \frac{|Q_j^k|^{mq[(1-\frac{\alpha}{mn})p_1^\prime-1](1-\frac{\alpha}{mn})}}{\nu_{\vec{\omega}}^q(Q_j^k)^{p_1^\prime(1-\frac{\alpha}{mn})-1} \prod_{i=1}^m \sigma_i(Q_j^k)^{\frac{qp_1^\prime}{p_i^\prime}(1-\frac{\alpha}{mn})}}\\
& \cdot \left( \prod_{i=1}^m \int _{Q_j^k} |f_i|\sigma_i \; dy_i\right)^q\\
\leq& C a^q [\vec{\omega}]_{A_{(\vec{P},q)}}^{p_1^\prime(1-\frac{\alpha}{mn})} \sum_{k,j} \frac{|E_j^k|^{mq[(1-\frac{\alpha}{mn})p_1^\prime-1](1-\frac{\alpha}{mn})}}{\nu_{\vec{\omega}}^q(Q_j^k)^{p_1^\prime(1-\frac{\alpha}{mn})-1} \prod_{i=1}^m \sigma_i(Q_j^k)^{\frac{qp_1^\prime}{p_i^\prime}(1-\frac{\alpha}{mn})}} \left(\prod_{i=1}^m \int _{Q_j^k} |f_i|\sigma_i \; dy_i\right)^q.\\
\end{align*}

We obverse that

$$
|E_j^k|=\int_{E_j^k} \nu_{\vec{\omega}}^{p\cdot\frac{1}{mp}} \sigma_1^{\frac{1}{mp_1^\prime}}\cdots \sigma_m^{\frac{1}{mp_m^\prime}} \leq  \nu_{\vec{\omega}}^p(E_j^k)^{\frac{1}{mp}} \prod_{i=1}^m \sigma_i(E_j^k)^{\frac{1}{mp_i^\prime}},
$$
where $\sigma_i=\omega_i^{-p_i^\prime}$.

Notice $p<q$ and use H\"{o}lder's inequality,

$$
|E_j^k| \leq \nu_{\vec{\omega}}^q(E_j^k)^{\frac{1}{mq}} |E_j^k|^{(\frac{1}{p}-\frac{1}{q})\frac{1}{m}} \prod_{i=1}^m \sigma_i(E_j^k)^{\frac{1}{mp_i^\prime}}.
$$
So,
$$
|E_j^k| \leq \nu_{\vec{\omega}}^q(E_j^k)^{\frac{1}{mq(1-\frac{\alpha}{mn})}} \prod_{i=1}^m \sigma_i(E_j^k)^{\frac{1}{mp_i^\prime(1-\frac{\alpha}{mn})}}.
$$
It's easy to see
$$
\frac{q(p_1^\prime-1)}{p_i^\prime} - \frac{q}{p_i} = \frac{qp_1^\prime}{p_i^\prime} - q \geq 0.
$$
Since $E_j^k\subset Q_j^k$, we have
$\sigma_i(E_j^k)^{\frac{q(p_1^\prime-1)}{p_i^\prime} - \frac{q}{p_i}}\leq \sigma_i(Q_j^k)^{\frac{qp_1^\prime}{p_i^\prime} - q}$, for $i=1,\cdots ,m$ and notice that
$\frac{m}{p_1}\geq \frac{1}{p}>\frac{\alpha}{n}$, so we have $p_1^\prime(1-\frac{\alpha}{mn})\geq 1$. Then
$\nu_{\vec{\omega}}^q(E_j^k)^{p_1^\prime(1-\frac{\alpha}{mn})-1} \leq \nu_{\vec{\omega}}^q(Q_j^k)^{p_1^\prime(1-\frac{\alpha}{mn})-1}$.
Therefore,
\begin{align*}
\int_{\R^n} &\mathcal{M}_\alpha^d(f_1 \sigma_1,\cdots ,f_m \sigma_m)(x)^q \nu_{\vec{\omega}}(x)^q \; dx\\
\leq& C [\vec{\omega}]_{A_{(\vec{P},q)}}^{p_1^\prime(1-\frac{\alpha}{mn})} \sum_{k,j} \prod_{i=1}^m \frac{\sigma_i(E_j^k)^{\frac{q}{p_i^\prime}[(1-\frac{\alpha}{mn})p_1^\prime-1]}} {\sigma_i(Q_j^k)^{\frac{qp_1^\prime}{p_i^\prime}(1-\frac{\alpha}{mn})}} \left( \prod_{i=1}^m \int _{Q_j^k} |f_i|\sigma_i \; dy_i\right)^q\\
\leq& C [\vec{\omega}]_{A_{(\vec{P},q)}}^{p_1^\prime(1-\frac{\alpha}{mn})} \sum_{k,j} \prod_{i=1}^m \left(\frac{1}{\sigma_i(Q_j^k)^{1-\frac{\alpha}{mn}}} \int_{Q_j^k} |f_i|\sigma_i \; dy_i \right)^q \sigma_i(E_j^k)^{\frac{q}{q_i}}\\
\leq& C [\vec{\omega}]_{A_{(\vec{P},q)}}^{p_1^\prime(1-\frac{\alpha}{mn})} \sum_{k,j} \prod_{i=1}^m \left(\int_{E_j^k} M_{\sigma_i,\alpha/m} (f_i) ^{q_i} \sigma_i \; dy_i \right)^{\frac{q}{q_i}}\\
\leq& C [\vec{\omega}]_{A_{(\vec{P},q)}}^{p_1^\prime(1-\frac{\alpha}{mn})} \prod_{i=1}^m \left( \sum_{k,j}\int_{E_j^k} M_{\sigma_i,\alpha/m} (f_i) ^{q_i} \sigma_i \; dy_i\right)^{\frac{q}{q_i}}\\
\leq& C [\vec{\omega}]_{A_{(\vec{P},q)}}^{p_1^\prime(1-\frac{\alpha}{mn})} \prod_{i=1}^m \left( \int_{\R^n} M_{\sigma_i,\alpha/m} (f_i) ^{q_i} \sigma_i \; dy_i\right)^{\frac{q}{q_i}}\\
\leq& C [\vec{\omega}]_{A_{(\vec{P},q)}}^{p_1^\prime(1-\frac{\alpha}{mn})} \prod_{i=1}^m ||f_i||_{L^{p_i}(\sigma_i)}^q,
\end{align*}
where $\frac{1}{q_i}=\frac{1}{p_i}-\frac{\alpha}{mn}$ and $M_{\sigma_i,\alpha/m} (f_i)(x) = \sup_{Q\ni x} \frac{1}{\sigma_i(Q)^{1-\frac{\alpha}{nm}}} \int_Q |f_i(y_i)|\sigma_i(y_i) \;dy_i$ and we use H\"{o}lder's inequality, and the boundedness of the weighted fractional maximal type operator.

Hence,
$$
||\mathcal{M}_\alpha(\vec{f})||_{L^q(\nu_{\vec{\omega}}^q)} \leq C_{n,m,\vec{P},q} [\vec{\omega}]_{A_{(\vec{P},q)}}^{\frac{p_1^\prime}{q}(1-\frac{\alpha}{mn})} \prod_{i=1}^m ||f_i ||_{L^{p_i}(\omega_i^{p_i})}
$$
Therefore, we obtain the upper bound of $\gamma$. Furthermore, the power is sharp. We complete the proof.\\

\noindent
\text{\it Proof of Theorem \ref{thm 1.2}.}

Applying (\ref{e 2.2}), we have
$$
||\mathcal{M}_{\Omega,\alpha}(\vec{f})||_{L^q(\nu_{\vec{\omega}}^q)}
\leq C_{n,m,\vec{P},q,s}  ||\Omega||_{L^s((S^{n-1})^m)}  ||\mathcal{M}_{\alpha s^\prime}(\vec{f}^{s^\prime})||_{L^{q/s^\prime}(\nu_{\vec{\omega}^{s^\prime}}^{q/s^\prime})}^{1/s^\prime}.
$$
According to Corollary \ref{1.4}, one obtains
$$
||\mathcal{M}_{\alpha s^\prime}(\vec{f}^{s^\prime})||_{L^{q/s^\prime}(\nu_{\vec{\omega}^{s^\prime}}^{q/s^\prime})}
\leq C_{n,m,\vec{P},q,s} [\vec{\omega^{s^\prime}}]_{A_{(\vec{P}/s^\prime,q/s^\prime)}}^{\gamma s^\prime}\prod_{i=1}^m ||f_i||_{L^{p_i}(\omega_i)}^{s^\prime}.
$$
where $\gamma$ is the same as the power in (\ref{e 1.7}), and satisfies\\
$(i)$if $1\leq s^\prime <p_1,\cdots ,p_m<\infty$, then $\frac{s^\prime}{q}\frac{mp/s^\prime}{mp/s^\prime-1}(1-\frac{\alpha s^\prime}{mn})\leq \gamma s^\prime \leq \frac{s^\prime}{q}(1-\frac{\alpha s^\prime}{mn})\max\{(\frac{p_1}{s^\prime})^\prime,\cdots ,(\frac{p_m}{s^\prime})^\prime\}$;\\
$(ii)$if $1\leq s^\prime<p_1=\cdots=p_m=s^\prime r<\infty$, then we have $s^\prime \gamma = \frac{r^\prime s^\prime}{q}(1-\frac{\alpha s^\prime}{mn})$.
Then, we complete the proof of Theorem \ref{thm 1.2}.\\

\noindent
\text{\it Proof of Theorem \ref{thm 1.5}.}

For any $\lambda
> 0$ and $k > 0$, denote that
\begin{center}$E_{\lambda} = \big\{x\in \mathbb{R}^n : \mathcal{M}_{\alpha}(\vec{f})(x) > \lambda\big\}$
 and $E_{\lambda, k} = E_{\lambda} \cap B(0,
k)$,\end{center}where $B(0,k) = \big\{x\in \mathbb{R}^n : |x|
\leqslant k\big\}$.

By the definition of $\mathcal{M}_\alpha$,
for any $x\in E_{\lambda, k}$ there is a cube $Q_x$
containing $x$ such that
\begin{equation}\label{e 3.1}\prod_{i = 1}^m\frac{1}{|Q_x|^{1 -
\frac{\alpha}{mn}}}\int_{Q_x}|f_i|\ > \lambda.\end{equation} Note that
$E_{\lambda, k} \subset \bigcup_{x\in E_{\lambda, k}}Q_{x}$. Then using
the Besicovitch covering lemma, we know that there exists an at most countable
subcollection of cubes ${Q_{x_j}}$ and a constant $C_n$ depending
only on the dimension such that
\begin{center}$E_{\lambda, k} \subset \displaystyle\bigcup_jQ_{x_j}$, \quad\quad
 $\displaystyle\sum_j\chi_{Q_{x_j}}(x) \leqslant C_n$.\end{center}
We summarize the proof into three subcases:\\
Case (i). If each $p_i = 1$, by the Definition \ref{def 1.3}, H\"{o}lder's inequality
and bounded overlap of $Q_{x_j}$, one obtain
\begin{eqnarray*}{\bigg(\int_{E_{\lambda, k}}{\nu_{\vec{\omega}}}^{\frac{n}{m n
- \alpha}}\,\bigg)}^{\frac{m n - \alpha}{mn}}
& \leqslant &
\sum_j{\bigg(\int_{Q_{x_j}}{\nu_{\vec{\omega}}}^{\frac{n}{m n -
\alpha}}\,\bigg)}^{\frac{m n -
\alpha}{mn}}{\bigg(\frac{1}{\lambda}\prod_{i =
1}^m\frac{(\inf\limits_{Q_{x_j}}w_i)^{-1}}{|Q_{x_j}|^{1 -
\frac{\alpha}{mn}}}\int_{Q_{x_j}}\big|f_iw_i\big|\,\bigg)}^{\frac{1}{m}}\\
&  \leqslant & C [\vec{w}]_{A_{(\vec {1},\frac{n}{mn-\alpha})}}^{\frac {mn-\alpha}{mn}}\sum_j {\bigg(\frac{1}{\lambda}\prod_{i =
1}^m\int_{Q_{x_j}}\big|f_i\omega_{i}\big|\,\bigg)}^{\frac{1}{m}}\\
& \leqslant & C[\vec{w}]_{A_{(\vec {1},\frac{n}{mn-\alpha})}}^{\frac {mn-\alpha}{mn}} {\bigg(\frac{1}{\lambda}\prod_{i =
1}^m\int_{\mathbb{R}^n}\big|f_i\omega_{i}\big|\,\bigg)}^{\frac{1}{m}}.\end{eqnarray*}

Case (ii). If each $p_i > 1$, by the fact that $\frac{p}{q} < 1$,
(\ref{e 3.1}) and H\"{o}lder's inequality and bounded overlap of
$Q_{x_j}$, we then
have\begin{eqnarray*}{\bigg(\int_{E_{\lambda,
k}}{\nu_{\vec{\omega}}}^q\bigg)}^{\frac{p}{q}}
& \leqslant &
\sum_j{\bigg(\int_{Q_{x_j}}{\nu_{\vec{\omega}}}^q\bigg)}^{\frac{p}{q}}{\bigg(\frac{1}{\lambda}\prod_{i
= 1}^m\frac{1}{|Q_{x_j}|^{1 -
\frac{\alpha}{mn}}}\int_{Q_{x_j}}\big|f_i\big|\bigg)}^{p}\\& \leqslant &
\sum_j{\bigg(\int_{Q_{x_j}}{\nu_{\vec{\omega}}}^q\bigg)}^{\frac{p}{q}}\frac{1}{\lambda^{p}}\prod_{i
= 1}^m\frac{1}{|Q_{x_j}|^{p - \frac{p
\alpha}{mn}}}{\bigg(\int_{Q_{x_j}}\big|f_i{\omega}_i\big|^{p_i}\bigg)}^{\frac{p}{p_i}}
{\bigg(\int_{Q_{x_j}}\omega_i^{-p'_i}\bigg)}^{\frac{p}{p'_i}}\\
& \leqslant & C\frac{[\vec{w}]_{A_{(\vec {P},{q})}}^{\frac {p}{q}}}{\lambda^{p}} \prod_{i =
1}^m{\bigg(\int_{\mathbb{R}^n}{\big|f_i{\omega}_i\big|}^{p_i}\bigg)}^{\frac{p}{p_i}}.\end{eqnarray*}

Case (iii) If there are some $p_i = 1$ and some $p_l > 1$, We assume $p_{l} > 1$ and $p_{i} = 1$. We combine the above methods in Case (i) and (ii) and obtain
that
\begin{eqnarray*}{\bigg(\int_{E_{\lambda,
k}}{\nu_{\vec{\omega}}}^q\bigg)}^{\frac{1}{q}}
& \leqslant &
\sum_j{\bigg(\int_{Q_{x_j}}{\nu_{\vec{\omega}}}^q\bigg)}^{\frac{1}{q}}\frac{1}{\lambda}\prod_{l}
\bigg\{\frac{1}{|Q_{x_j}|^{1 -
\frac{\alpha}{mn}}}{\Big(\int_{Q_{x_j}}\big|f_{l}{\omega}_{l}\big|^{p_{l}}\Big)}^{\frac{1}{p_{l}}}
{\Big(\int_{Q_{x_j}}\omega_{l}^{-p'_{l}}\Big)}^{\frac{1}{p'_{l}}}\bigg\}\\
& & \times\prod_{i} {\bigg\{\frac{1}{|Q_{x_j}|^{1 -
\frac{\alpha}{mn}}}{\Big(\inf_{Q_{x_j}}\omega_{i}\Big)}
^{-1}\int_{Q_{x_j}}\big|f_{i}\omega_{i}\big|\bigg\}}\\
& \leqslant & \frac{C[\vec{w}]_{A_{(\vec {P},{q})}}^{\frac {1}{q}}}{\lambda} \prod_{i =
1}^m{\bigg(\int_{\mathbb{R}^n}{\big|f_i{\omega}_i\big|}^{p_i}\bigg)}^{\frac{1}{p_i}}.
\end{eqnarray*}

It is obvious that all $C$ are independent of $k$, which implies
(\ref{e 1.15}) by
monotone convergence theorem.

\section{The proof of Theorem \ref{thm 1.6}}
\setcounter{equation}{0}

For a general Banach function space $X$, we know that
$$
\int _{\R^n} |f(x) g(x)|\; dx \leq ||f||_X ||g||_{X^\prime},
$$
which we will use in the proof of Theorem \ref{thm 1.6}, where $X^\prime$ consists of measurable functions $f$ for which satisfies
$$
||f||_{X^\prime}=\sup_{||g||_X\leq 1} \int_{\R^n} |f(x) g(x)|\; dx <\infty.
$$

\noindent
{\it Proof of Theorem \ref{thm 1.6}.}
According to the proof of Corollary \ref{1.3}, we only need to show $\mathcal{M}_\alpha^d$ instead of $\mathcal{M}_\alpha$. Let $\Omega_k=\{x\in \R^n :\mathcal{M}_\alpha^d(\vec{f})(x)>a^k\}=\cup_j Q_j^k$, where $a=2^{m(n+1)}$. Then
Note that $p<q$, and by the H\"{o}lder's inequality, we have
\begin{align*}
\int_{\R^n} &\mathcal{M}_\alpha^d(x)^q u(x) \; dx
\leq a^q \sum_{k,j} \left(\prod_{i=1}^m \frac{1}{|Q_j^k|^{1-\frac{\alpha}{mn}}} \int _{Q_j^k} |f_i v_i v_i^{-1}| \; dy_i\right)^q u(Q_j^k)\\
\leq& a^q \sum_{k,j} \left(\prod_{i=1}^m \frac{1}{|Q_j^k|} \int _{Q_j^k} |f_i v_i v_i^{-1}| \; dy_i\right)^q \frac{u(Q_j^k)}{|Q_j^k|} |Q_j^k|^{\frac{q}{p}}\\
\leq& a^q \sum_{k,j} \left(|Q_j^k|\prod_{i=1}^m ||f_i v_i||^p_{X_i^\prime,Q_j^k}\right)^{\frac{q}{p}} \prod_{i=1}^m ||v_i^{-1}||^q_{X_i,Q_j^k} \frac{u(Q_j^k)}{|Q_j^k|} \\
\leq& 2^{\frac{q}{p}} a^q K^q \sum_{k,j} \left(\int_{E_j^k} \prod_{i=1}^m M_{X_i^\prime}(f_i v_i)(x)^p \; dx \right)^{\frac{q}{p}}
\leq 2^{\frac{q}{p}} a^q K^q \left(\sum_{k,j} \int_{E_j^k} \prod_{i=1}^m M_{X_i^\prime}(f_i v_i)(x)^p \; dx \right)^{\frac{q}{p}}\\
\leq& C ^q K^q \prod_{i=1}^m \left(\int_{\R^n} M_{X_i^\prime}(f_i v_i)(x)^{p_i} \; dx \right)^{\frac{q}{p_i}}
\leq C ^q K^q \prod_{i=1}^m ||M_{X_i^\prime}||_{p_i}^q ||f_i v_i||_{p_i}^q,
\end{align*}
So we obtain that
$$
||\mathcal{M}_\alpha(\vec{f})||_{L^q(u)} \leq CK \prod _{i=1}^m ||M_{X_i^\prime}||_{L^{p_i}(\R^n)} ||f_i v_i||_{L^{p_i}(\R^n)}.
$$
The proof is ended.

\remark \label{rem 4.1}
The result of Theorem \ref{thm 1.6} can be seen as a two weights version of (\ref{e 1.8}).

\section{The proof of Corollary {1.7} and Theorem \ref{thm 1.8}-\ref{thm 1.9}}
\setcounter{equation}{0}

\begin{lemma}[Pointwise estimate of $\mathcal{I}_{\Omega, \alpha}$]Suppose that $\epsilon > 0$
satisfying $0 < \alpha - \epsilon$. Then for any $x\in
\mathbb{R}^n$, there is a constant $C > 0$ independent of $\vec{f}$
such that
\begin{equation}\label{pointwise kernel}\big|\mathcal{I}_{\Omega,
\alpha}(\vec{f})(x)\big| \leqslant C {\big(\mathcal{M}_{\Omega,
\alpha +
\epsilon}(\vec{f})(x)\big)}^{\frac{1}{2}}{\big(\mathcal{M}_{\Omega,
\alpha -
\epsilon}(\vec{f})(x)\big)}^{\frac{1}{2}}.\end{equation}
\end{lemma}
\textbf{Proof.} For any cube $Q\ni
x$, we have$$\aligned\mathcal {I}_{\Omega, \alpha}(\vec{f})(x)
&= \displaystyle \int_{Q^m}\frac{|\Omega(x-y_1,\cdots x - y_m)|\prod_{i = 1}^m|f_i(y_i)|}{{|(x -
y_1, \cdots, x - y_m)|}^{mn - \alpha}}\,d\vec{y}\\&\quad+
\int_{{(\mathbb{R}^n)}^m \backslash Q^m}\frac{|\Omega(x-y_1,\cdots x - y_m)|\prod_{i =
1}^m|f_i(y_i)|}{{|(x - y_1, \cdots, x - y_m)|}^{mn -
\alpha}}\,d\vec{y}\\&\triangleq J_1 + J_2,\endaligned
$$

where $J_1 =
\int_{Q^m}\frac{|\Omega(x-y_1,\cdots x - y_m)|\prod_{i = 1}^m|f_i(y_i)|}{{|(x -
y_1, \cdots, x - y_m)|}^{mn - \alpha}}\,d\vec{y}$, $J_2 =
\int_{{(\mathbb{R}^n)}^m \backslash Q^m}\frac{|\Omega(x-y_1,\cdots x - y_m)|\prod_{i =
1}^m|f_i(y_i)|}{{|(x - y_1, \cdots, x - y_m)|}^{mn -
\alpha}}\,d\vec{y}$.

First, for $J_1$, we have \begin{eqnarray*}J_1 & \leqslant & \displaystyle C\sum_{j = 0}^\infty{{(2^{- j
-1}{|Q|}^{\frac{1}{n}})}^{-(mn - \alpha)}}\int_{{(2^{- j}Q)}^m \backslash {(2^{- j -
1}Q)}^m}|\Omega(x - y_1,\cdots,x-y_m)|\prod_{i
= 1}^m|f_i(y_i)|\,d\vec{y}\\
& \leqslant & C {|Q|}^{\frac{\epsilon}{n}}\mathcal{M}_{\Omega,
\alpha - \epsilon}(\vec{f})(x).\end{eqnarray*}

For $J_2$, similarly, we have$$\aligned
J_2 &\le \displaystyle\sum_{j =
0}^\infty\int_{{(2^{j}Q)}^m \backslash {(2^{j -
1}Q)}^m}\frac{|\Omega(x-y_1,\cdots x - y_m)|\prod_{i =
1}^m|f_i(y_i)|}{{|(x - y_1, \cdots, x - y_m)|}^{mn -
\alpha}}\,d\vec{y}\\& \le C {|Q|}^{\frac{- \epsilon}{n}}\mathcal{M}_{\Omega,
\alpha + \epsilon}(\vec{f})(x).
\endaligned
$$

By taking ${|Q|}^{\frac{2\epsilon}{n}} = \frac{\mathcal{M}_{\Omega,
\alpha + \epsilon}(\vec{f})(x)}{\mathcal{M}_{\Omega, \alpha -
\epsilon}(\vec{f})(x)}$, we have (\ref{pointwise kernel}).\\

\noindent
{\it Proof of Corollary \ref{cor 1.7}.} Note that
$\frac{q}{2q_{\epsilon}} + \frac{q}{2q_{- \epsilon}} = 1$, then by
(\ref{pointwise kernel}), the H\"{o}lder's inequality, we
have\begin{eqnarray*}\lefteqn{{\bigg(\int_{\mathbb{R}^n}{\Big(\mathcal{I}_{\Omega,
\alpha}(\vec{f})\big|
\nu_{\vec{\omega}}\Big)}^{q}\,\bigg)}^{\frac{1}{q}}}\\
&  \leqslant & C {\bigg(\int_{\mathbb{R}^n}\Big({\mathcal{M}_{\Omega,
\alpha +
\epsilon}(\vec{f})\nu_{\vec{\omega}}\Big)}^{q_{\epsilon}}\,\bigg)}^{\frac{1}{2q_{\epsilon}}}
{\bigg(\int_{\mathbb{R}^n}\Big({\mathcal{M}_{\Omega, \alpha -
\epsilon}(\vec{f})\nu_{\vec{\omega}}\Big)}^{q_{-
\epsilon}}\,\bigg)}^{\frac{1}{2q_{- \epsilon}}}.
\end{eqnarray*}
Then Corollary \ref{cor 1.7} follows from (\ref{e 1.6}) and (\ref{e 1.9}).\\

\noindent
{\it Proof of Theorem \ref{thm 1.8}.}
For any cube $Q\subset \R^n$, suppose $B$ with radius $r$ is the circumscribed ball of the cube $Q$ and denote that $\vec{f}^{s^\prime}=(f_1^{s^\prime},\cdots ,f_m^{s^\prime})$. Then we have the following estimate

\begin{align}
\mathcal{M}_\Omega(\vec{f})(x) &\leq \sup_{Q \ni 0}  \left( \frac{1}{|Q|^m} \int_{Q^m} \prod_{i=1}^m|f_i(x-y_i)|^{s^\prime} \; d\vec{y}  \right)^{1/s^\prime}  \left( \frac{1}{|Q|^m} \int_{B^m} |\Omega(y_1,\cdots ,y_m)|^s\; d\vec{y} \right)^{1/s}\nonumber\\
&\leq C ||\Omega||_{L^s((S^{n-1})^m)} \sup_{Q \ni 0} \prod_{i=1}^m  \left(\frac{1}{|Q|} \int_Q |f_i(x-y_i)|^{s^\prime} \; dy_i \right)^{1/s^\prime}\nonumber\\
&= C ||\Omega||_{L^s((S^{n-1})^m)}  \mathcal{M}(\vec{f}^{s^\prime})(x)^{\frac{1}{s^\prime}}.\label{e 5.2}
\end{align}

Apply Theorem A, and notice $\vec{\omega} \in A_{\vec{P}/s^\prime}$, we can see

\begin{align*}
||\mathcal{M}_\Omega(\vec{f})||_{L^p(\nu_{\vec{\omega}})}
\leq& C_{n,m,\vec{P},s}  ||\Omega||_{L^s((S^{n-1})^m)}  ||\mathcal{M}(\vec{f}^{s^\prime})||_{L^{p/s^\prime}(\nu_{\vec{\omega}})}^{1/s^\prime}\\
\leq& C_{n,m,\vec{P},s}  ||\Omega||_{L^s((S^{n-1})^m)}   [\vec{\omega}]_{A_{\vec{P}/s^\prime}}^{1/p} \prod_{i=1}^m [\omega_i^{1-(\frac{p_i}{s^\prime})^\prime}]_{A_\infty}^{1/p_i} \prod_{i=1}^m ||f_i||_{L^{p_i}(\omega_i)}.
\end{align*}
So we complete the proof of Theorem \ref{thm 1.8}.
The proof of Theorem \ref{thm 1.9} is more easy, now we will give it.
By (\ref{e 5.2}) and applying Theorem B, we know

\begin{align*}
||\mathcal{M}_\Omega(\vec{f})||_{L^p(\nu_{\vec{\omega}})}
\leq& C_{n,m,\vec{P},s}  ||\Omega||_{L^s((S^{n-1})^m)} ||\mathcal{M}(\vec{f}^{s^\prime})||_{L^{p/s^\prime}(\nu_{\vec{\omega}})}^{1/s^\prime}\\
\leq& C_{n,m,\vec{P},s}  ||\Omega||_{L^s((S^{n-1})^m)} [\vec{\omega}]_{A_{\vec{P}/s^\prime}}^\gamma \prod_{i=1}^m ||f_i||_{L^{p_i}(\omega_i)},
\end{align*}
where $\gamma$ satisfied
$$
\frac{m}{mp/s^\prime-1} \frac{1}{s^\prime} \leq \gamma \leq \frac{1}{s^\prime} \max\{\frac{(p_1/s^\prime)^\prime}{p/s^\prime},\cdots ,\frac{(p_m/s^\prime)^\prime}{p/s^\prime}\}.
$$
That is,

$$
\frac{m}{mp-s^\prime} \leq \gamma \leq \max\{\frac{(p_1/s^\prime)^\prime}{p},\cdots ,\frac{(p_m/s^\prime)^\prime}{p}\}.
$$
We complete the proof of the theorem.

\end{document}